\documentclass[12pt,psamsfonts]{amsart}

\usepackage{amsmath,amssymb}
\usepackage{graphicx}
\usepackage{color}
\usepackage{psfrag} % See the section on figures.
\usepackage{overpic} % See the section on figures.
\usepackage{array} % Function notation
 
\newtheorem{theorem}{Theorem}
\newtheorem{proposition}[theorem]{Proposition}

\newtheorem{lemma}[theorem]{Lemma}
\newtheorem {remark}[theorem]{Remark}
\newtheorem{definition}[theorem]{Definition}

\title[Nilpotent singular points with Odd Andreev number]{Monodromic nilpotent singular points with odd Andreev number and the Center problem}
\author[L. Queiroz and C. Pessoa]{}


\begin{document}
 \maketitle

\centerline{\scshape  Claudio Pessoa,  \; Lucas Queiroz}
\medskip

{\footnotesize \centerline{Universidade Estadual Paulista (UNESP), Instituto de Bioci\^encias Letras e Ci\^encias Exatas,} \centerline{R. Cristov\~ao Colombo, 2265, 15.054-000, S. J. Rio Preto, SP, Brasil }
\centerline{\email{c.pessoa@unesp.br} and \email{lucas.queiroz@unesp.br}}}

\medskip

\bigskip

\begin{quote}{\normalfont\fontsize{8}{10}\selectfont
{\bfseries Abstract.}
Given a nilpotent singular point of a planar vector field, its monodromy is associated with its Andreev number $n$. The parity of $n$ determines whether the existence of an inverse integrating factor implies that the singular point is a nilpotent center. For $n$ odd, this is not always true. We give a characterization for a family of systems having Andreev number $n$ such that the center problem cannot be solved by the inverse integrating factor method. Moreover, we study general properties of this family, determining necessary center conditions for every $n$ and solving the center problem in the case $n=3$.
\par}
\end{quote}

\section{Introduction and Preliminary results}
%\label{sec:01}

The \emph{Center Problem} or \emph{Center-focus Problem} for an analytical planar vector field consists of distinguishing whether a monodromic singular point is either a center or a focus. The study of this problem started with the works of Poincaré and Dulac in the last century, but it remains open until the present day and several other related problems have arisen from it, e.g. the isochronicity and the cyclicity problems (see \cite{Romanovski}). Here we are interested in the version of this problem for nilpotent monodromic singular points.

An analytical planar vector field has a nilpotent singular point if its Jacobian matrix, at the singular point, have all null eigenvalues but is not the null matrix. By means of a linear change of variables, the system associated to a vector field with this property can be written in the form:
\begin{equation}\label{FNeq1}
	\begin{array}{lcr}
		\dot{x}=y+P(x,y),\\
		\dot{y}=Q(x,y),
	\end{array}
\end{equation}
where $P, Q$ are analytic and $j^1P(0)=j^1Q(0)=0$. Here, $j^lF(0)$ denotes the $l$-jet of a smooth function $F$ about the origin. In \cite{Moussu}, the authors establish a characterization of nilpotent centers via time-reversibility (see also, \cite{Zoladek}). However, the center problem for nilpotent singular points is far from being solved, since there are only few families for which the center conditions are known (see, for instance, \cite{Algaba,AlvarezGasull1,AlvarezGasull2}).

For nilpotent singular case, in order to study the center problem, we first need to identify which vector fields having a nilpotent singular point are monodromic. Andreev's result \cite{Andreev} gives us this first characterization.

\begin{theorem}[Andreev's Theorem]\label{TeoAndreev}
	Let $X$ be the vector field associated to system \eqref{FNeq1} and the origin be an isolated singular point. Let $y=F(x)$ be the solution of the equation $y+P(x,y)=0$ in a neighborhood of $(0,0)$ and consider $f(x)=Q(x,F(x))$ and $\Phi(x)=\mbox{\rm div}X\vert_{(x,F(x))}$. Then, we can write
	$$f(x)=ax^{\alpha}+O(x^{\alpha+1}),$$
	$$\Phi(x)=bx^{\beta}+O(x^{\beta+1}).$$
	The origin is monodromic if and only if $a<0,\alpha=2n-1$ and one of the following conditions holds:
	\begin{itemize}
		\item[i)]$\beta>n-1$;
		\item[ii)]$\beta=n-1$ and $b^2+4an<0$;
		\item[iii)]$\Phi\equiv 0$.
	\end{itemize}
\end{theorem}

So we define the \emph{Andreev number} of a nilpotent singular point by the number $n$ in function $f(x)=ax^{2n-1}+O(x^{2n})$ above.

In \cite{GarciaFiI} the author proves that the Andreev number is invariant by analytical orbital equivalence, i.e. by diffeomorphisms and time rescalings. In the same paper, a method to solve the center problem for nilpotent systems is presented which consists in finding a formal inverse integrating factor for such systems, i.e. a formal series $V(x,y)$ satisfying $XV=V\mbox{div}X$. Throughout this work, we will refer to this method as the \emph{inverse integrating factor method}. More precisely, the method is based on the following results:

\begin{theorem}[Theorem 3 in \cite{GarciaFiI}]%\label{TeoGarciaFiI}
	Consider system \eqref{FNeq1} with a monodromic nilpotent singular point at the origin such that the Andreev number $n$ is even. If \eqref{FNeq1} has a formal inverse integrating factor then the origin is a center. Moreover if $n=2$ then the converse is also true. 
\end{theorem}

The case with Andreev number $n=2$ was also studied in \cite{LiuLi3Ord,LiuLiNew1} where the authors also found a canonical form for nilpotent monodromic singularities of third order. The existence of a formal inverse integrating factor by itself is not sufficient to prove that the singular point is a nilpotent center. The best result in this direction, also presented in \cite{GarciaFiI}, is the next theorem. Before stating it, we need the following definition:

\begin{definition}%\label{defQuasihomogeneous}
A polynomial $p\in\mathbb{R}[x,y]$ is a \emph{$(t_1,t_2)$-quasi-\-homogeneous polynomial of weighted degree $k$} if $p(\lambda^{t_1} x, \lambda^{t_2} y)=\lambda^kp(x,y)$. A general expression for such polynomials is $p(x,y)=\sum_{t_1i+t_2j=k}a_{ij}x^iy^j$ where $a_{ij}\in\mathbb{R}$. The vector space of all $(t_1,t_2)$-quasi-\-homogeneous polynomial of weighted degree $k$ is denoted by $\mathcal{P}^{(t_1,t_2)}_{k}$.
\end{definition}

\begin{theorem}[Theorem 4 in \cite{GarciaFiI}]\label{TeoGarciaFiI2}
	Consider system \eqref{FNeq1} with a monodromic nilpotent singular point at the origin and Andreev number $n$ satisfying the condition $\beta>n-1$ given in Theorem \ref{TeoAndreev}. If there is a formal inverse integrating factor $V(x,y)=\sum_{j\geqslant 2n}V_j(x,y)$ where $V_j$ are $(1,n)$-quasi-homogeneous polynomials of weighted degree $j$ and $V_{2n}\neq 0$, then the origin is a center. The converse is not true.
\end{theorem}

Our main objective is to study families of polynomial systems of the form \eqref{FNeq1} having a singular point with odd Andreev number for which the inverse integrating factor method is not sufficient to solve the center problem. We characterize some of those families, exhibit some of their properties focusing on the center problem and obtain some general necessary center conditions.

\section{Statement of the main results}

When studying the nilpotent singular points of analytic vector fields, one will eventually stumble upon the definition of quasi-homogeneous polynomials. For systems having a singular point with Andreev number $n$, the $(1,n)$-quasi-\-homogeneous polynomials are especially important. Therefore, some definitions are in order:

\begin{definition}
	Consider an analytical mapping $V:\mathbb{R}^m\to\mathbb{R}^l$. We can write $V=\sum_{j\geqslant 0}V_j$ where all the components of $V_j$ are $(t_1,t_2)$-quasi-\-homogeneous polynomials of weighted degree $j$. The \emph{$k$-jet of $V$ by $(t_1,t_2)$-quasi-\-homogeneous terms} is given by the expression $\sum_{j=0}^{k}V_j$.
\end{definition}
\begin{definition}
	A vector field is a \emph{$(1,n)$-quasi-homogeneous polynomial vector field of maximum weighted degree $k$} if its components are sums of $(1,n)$-quasi-homogeneous polynomials of weighted degree less or equal to $k$.
\end{definition}

Expanding the functions $P,Q$ in system \eqref{FNeq1} in their power series about the origin, we obtain the following form:
\begin{equation}\label{FNeq2}
	\begin{array}{lcr}
		\dot{x}=y+\displaystyle\sum_{i+j\geqslant 2}a_{ij}x^iy^j,\\
		\dot{y}=\displaystyle\sum_{i+j\geqslant 2}b_{ij}x^iy^j,\\
	\end{array}
\end{equation} 
where $a_{ij},b_{ij}\in\mathbb{R}$. In the case which the origin is a monodromic nilpotent singular point of \eqref{FNeq2}, it is possible, via analytic changes of variables, to write \eqref{FNeq2} in the following form: 
\begin{equation}\label{FNfamilia}
	\begin{array}{lcr}
		\dot{x}=y+\mu x^n+\displaystyle\sum_{i+nj\geqslant n+1}\tilde{a}_{ij}x^iy^j,\\
		\dot{y}=-nx^{2n-1}+n\mu x^{n-1}y+\displaystyle\sum_{i+nj\geqslant 2n}\tilde{b}_{ij}x^iy^j,\\
	\end{array}
\end{equation}
where $\mu\in\mathbb{R}$. This proof of this fact is based on the normal preforms presented in \cite{AlgabaInte} and will be detailed in Section \ref{SecFamiliaOrdemn}. Note that, by means of the change of variables $u=\sqrt[2n-2]{n}\,x,\;v=-\sqrt[2n-2]{n}\,y$, the above system can be written as
\begin{equation}\label{FNfamilia2}
	\begin{array}{lcr}
		\dot{x}=-y+\tilde{\mu} x^n+\displaystyle\sum_{i+nj\geqslant n+1}\hat{a}_{ij}x^iy^j,\\
		\dot{y}=x^{2n-1}+n\tilde{\mu} x^{n-1}y+\displaystyle\sum_{i+nj\geqslant 2n}\hat{b}_{ij}x^iy^j.\\
	\end{array}
\end{equation}

The canonical form \eqref{FNfamilia2} is in the class of systems studied in paper \cite{GasullTorre} where the authors provide an easier method to compute focal values of any order by a recurrent formula. Moreover, this canonical form of monodromic systems having a singular point with odd Andreev number allow us to detect the nilpotent focus at the origin for $\mu\neq 0$ regardless the particular value of $n$. Since $\mu\neq 0$ holds if and only if the considered original systems satisfy monodromy condition (ii) in Theorem \ref{TeoAndreev}, we are able to prove the following result.

\begin{theorem}\label{TeoFocoNilpotenteForte}
	Suppose that system \eqref{FNeq1} has nilpotent monodromic singular point with odd Andreev number $n$. If the system satisfies the monodromy condition $\beta=n-1$ then the origin is a nilpotent focus.
\end{theorem}

%\begin{teo}\label{TeoFocoNilpotente2}
%	Consider system \eqref{FNfamilia} with odd Andreev number $n$, if $\mu\neq 0$ then the origin is a nilpotent focus.
%\end{teo}

The proof of this theorem will be presented in Section \ref{secPolar}. Theorem \ref{TeoFocoNilpotenteForte} does not hold for systems having even Andreev number. For instance, it is possible for system \eqref{FNfamilia} with $n=2$ and $\mu\neq 0$ to have a nilpotent center at the origin (see \cite{LiuLi3Ord}). Those types of singular points are called \emph{third order $0$-class critical points} and are studied in \cite{LiuLi3Ord,LiuLiNew1,GarciaFiI}. 

To further investigate monodromic singularities having odd Andreev number, we consider a particular family of systems \eqref{FNfamilia}. More precisely, we consider all $(1,n)$-quasi-homogeneous systems of maximum weighted degree $2n-1$:
\begin{equation*}
	\begin{array}{lcr}
		\dot{x}=y+\mu x^n+y\displaystyle\sum_{k=1}^{n-1}a_{k1}x^{k}+\displaystyle\sum_{k=n+1}^{2n-1}a_{k0}x^{k},\\
		\dot{y}=-nx^{2n-1}+n\mu x^{n-1}y.\\
	\end{array}
\end{equation*}
From Theorem \ref{TeoFocoNilpotenteForte} we know that for $\mu\neq 0$, the origin is a nilpotent focus. Thus, to study the center problem, it is more convenient to study the following system which is obtained by making the change of variables $u=\sqrt[2n-2]{n}\,x,\;v=-\sqrt[2n-2]{n}\,y$ in the above system for $\mu=0$.
\begin{equation}\label{FNeq1nhomo2n-1}
	\begin{array}{lcr}
		\dot{x}=-y+y\displaystyle\sum_{k=1}^{n-1}a_{k1}x^{k}+\displaystyle\sum_{k=n+1}^{2n-1}a_{k0}x^{k},\\
		\dot{y}=x^{2n-1}.\\
	\end{array}
\end{equation}

The next theorem will be proven in section \ref{secImpossibilidade}.

\begin{theorem}\label{TeoFNimpossibilidadeGarcia}
	Consider system \eqref{FNeq1nhomo2n-1}, with $n$ odd. If it admits formal inverse integrating factor $V(x,y)$, then either $a_{k0}=0$ for $k=n+1,\dots,2n-1$ or the $(1,n)$-quasi-homogeneous part of weighted degree $2n$ of $V(x,y)$ is null. 
\end{theorem}
From this result we conclude that the formal inverse integrating factor method applied to systems \eqref{FNeq1nhomo2n-1} with $n$ odd is not useful to determine center conditions. Nevertheless, the properties of such systems allow us to obtain general necessary center conditions.

\begin{proposition}\label{PropoFNcondcentro}
	Consider system \eqref{FNeq1nhomo2n-1} with $n$ odd. If the origin is a nilpotent center, then we must have $a_{n+2,0}=-a_{11}a_{n+1,0}$.
\end{proposition}

We give a proof of this result in section \ref{secImpossibilidade}. For $(1,3)$-quasi-\-homo\-geneous systems \eqref{FNeq1nhomo2n-1} of maximum weighted degree $5$ having a mono\-dromic nilpotent singular point with Andreev number $n=3$, we were able to completely solve the center problem which solution will be presented in section \ref{sec135}.

\begin{theorem}\label{TeoFNCentro135}
	Considering the following system having a nilpotent monodromic singular point at the origin with Andreev number $3$:
	\begin{equation}\label{FNeq13homo5}
		\begin{array}{lcr}
			\dot{x}=-y+\mu x^3+a_{11}xy+a_{21}x^2y+a_{40}x^4+a_{50}x^5,\\
			\dot{y}=x^{5}+3\mu x^{2}y.\\
		\end{array}
	\end{equation}
	The origin is a nilpotent center if and only if $\mu=a_{11}a_{40}=a_{50}=0$.
\end{theorem}

We will now develop the tools and results necessary to prove these theorems. The methods and studies highlighted in this paper can be used in a more broad context, and we believe they are useful in general.

\section{Classification of monodromic nilpotent singularities}\label{SecFamiliaOrdemn}

In this section we prove that if the origin of system \eqref{FNeq1} is an isolated monodromic nilpotent singular point, then it can be transformed into system \eqref{FNfamilia} via a local change of variables.

First we state the next lemma which provides a shortcut to verify monodromy under some useful hypothesis.

\begin{lemma}[Monodromy criterion]\label{LemaMonodromiaDelta}
	Consider system \eqref{FNeq1}, the functions $F(x),f(x),\Phi(x)$ given in Theorem \ref{TeoAndreev}. Let $\tilde{a}$ be the coefficient of $x^{2n-1}$ in the power series expansion of $f(x)$, $\tilde{b}$ the coefficient of $x^{n-1}$ in the power series expansion of $\Phi(x)$ and $\Delta=\tilde{b}^2+4\tilde{a}n$. Suppose $j^{2n-2}f(0)=0$ and $j^{n-2}\Phi(0)=0$, then the singular point is monodromic with Andreev number $n$ if and only if $\Delta<0$. 
\end{lemma}
\noindent\textbf{Proof: }Note that $\Delta<0$ gives us immediately the first monodromy condition, i.e. $\alpha=2n-1$ and $a=\tilde{a}<0$. If $\tilde{b}=0$, then either $j^r\Phi(0)\equiv 0$ or $\Phi(x)\in O(x^n)$ and the monodromy condition (i) in Theorem \ref{TeoAndreev} is satisfied. If $\tilde{b}\neq 0$, the condition (ii) in Theorem \ref{TeoAndreev} is satisfied. Thus, the singular point is monodromic with Andreev number $n$. Conversely, suppose the singular point is monodromic with Andreev number $n$. Then $\tilde{a}<0$ and one of the monodromy conditions (i) or (ii) in Theorem \ref{TeoAndreev} must hold. If (i) holds, $\tilde{b}$ is null and $\tilde{a}<0$ implies $\Delta<0$. On the other hand, if (ii) is satisfied then $\Delta$ must be negative. Therefore $\Delta<0$ regardless the monodromy condition.\qed 

\begin{proposition}\label{PropoPreformMono}
Consider system \eqref{FNeq1}. The origin is a monodromic isolated singular point if and only if there exists a local analytical change of variables that transforms \eqref{FNeq1} into \eqref{FNfamilia}.
\end{proposition}
\noindent\textbf{Proof: }Note that the origin of system \eqref{FNfamilia} is monodromic by Theorem \ref{TeoAndreev}. Conversely, considering the functions defined in Theorem \ref{TeoAndreev}, the analytic change of variables $x\to x$, $y\to y-F(x)$, transforms system \eqref{FNeq1} into
\begin{equation}\label{FNeqAndreev}
	\begin{array}{lcr}
		\dot{x}=y+y\tilde{P}(x,y),\\
		\dot{y}=f(x)+y\Phi(x)+y^2\tilde{Q}(x,y).\\
	\end{array}
\end{equation}
By Theorem \ref{TeoAndreev}, since the origin is monodromic, the above system can be rewritten as
\begin{equation*}
	\begin{array}{lcr}
		\dot{x}=y+\sum_{r\geqslant n+1}P_r(x,y),\\
		\bigskip
		\dot{y}=ax^{2n-1}+\tilde{b}x^{n-1}y+\sum_{r\geqslant 2n}Q_r(x,y),\\
	\end{array}
\end{equation*}
where $P_r(x,y),Q_r(x,y)\in\mathcal{P}^{(1,n)}_{r}$ and $a\neq 0$, $\tilde{b}\in\mathbb{R}$. Note that the origin has Andreev number $n$. By Lemma \ref{LemaMonodromiaDelta}, we have $\Delta=\tilde{b}^2+4an<0$. Define $A=\frac{\Delta}{4n^2}$. The mapping $$\Psi(x,y)=\left(|A|^{\frac{1}{2(n-1)}}x, |A|^{\frac{1}{2(n-1)}}\left(y-\frac{\tilde{b}}{2n}x^{n}\right)\right),$$
is a diffeomorphism, and it transforms \eqref{FNeqAndreev} into
\begin{equation*}
	\begin{array}{lcr}
		\dot{x}=y+\mu x^n+\sum_{r\geqslant n+1}\tilde{P}_r(x,y),\\
		\bigskip
		\dot{y}=-nx^{2n-1}+n\mu x^{n-1}y+\sum_{r\geqslant 2n}\tilde{Q}_r(x,y),\\
	\end{array}
\end{equation*}
where $\tilde{P}_r(x,y),\tilde{Q}_r(x,y)\in\mathcal{P}^{(1,n)}_{r}$ and $\mu=\frac{\tilde{b}}{2n|A|^{\frac{1}{2}}}$.\qed

\begin{remark}\label{obsbeta=n-1}
The expression of $\mu$ given above yields the following implication: Consider system \eqref{FNeq1} with a monodromic nilpotent singular point. It satisfies the monodromy condition $\beta=n-1$ (item (ii) from Theorem \ref{TeoAndreev}) if and only if $\mu\neq0$ in canonical form \eqref{FNfamilia}. 
\end{remark}

\section{Generalized polar coordinates}\label{secPolar}

We now turn our attention to systems in family \eqref{FNfamilia}. It is convenient to make the following change of variables, first proposed by Lyapunov \cite{Lyapunov}, $x=r\,\mbox{Cs}\,\theta,\;y=r^n\,\mbox{Sn}\,\theta$, where Cs $\theta$ and Sn $\theta$ are the solutions of the Cauchy problem:
$$\begin{array}{lcr}
	\dfrac{dx}{d\theta}=-y;\\
	\dfrac{dy}{d\theta}=x^{2n-1},\\
	x(0)=1, y(0)=0.
\end{array}$$

We refer to this change of variables as the \emph{generalized polar coordinates}. The functions Cs $\theta$ and Sn $\theta$ have properties analogous to the classic trigonometric functions, which make its algebraic manipulation easier. We state some of these properties whose proofs the reader can find in \cite{Lyapunov, GasullTorre}.

\begin{proposition}\label{propoPropriedadesTrigo}
	The following holds:
	\begin{itemize}
		\item[a)] {\rm Cs}\,$\theta$ and {\rm Sn}\,$\theta$ are $T$-periodic functions, where 
		$$T=2\sqrt{\dfrac{\pi}{n}}\dfrac{\Gamma(\frac{1}{2n})}{\Gamma(\frac{n+1}{2n})};$$
		\item[b)] $\mbox{\rm Cs}^{2n}\theta+n\mbox{\rm Sn}^2\theta=1$;
		\item[c)] {\rm Cs}\,$\theta$ is an even function and {\rm Sn}\,$\theta$ is an odd function;
		\item[d)] {\rm Cs}\,$(\frac{T}{2}+\theta)=-\mbox{\rm Cs}\,\theta$ and {\rm Sn}\,$(\frac{T}{2}+\theta)=-\mbox{\rm Sn}\,\theta$;
		\item[e)] {\rm Cs}\,$(\frac{T}{2}-\theta)=-\mbox{\rm Cs}\,\theta$ and {\rm Sn}\,$(\frac{T}{2}-\theta)=\mbox{\rm Sn}\,\theta$;
		\item[f)] $\displaystyle\int_{0}^{\theta}\mbox{\rm Sn}\,\varphi\;\mbox{\rm Cs}^q\varphi \;d\varphi=-\dfrac{\mbox{\rm Cs}^{q+1}\theta}{q+1}$;
		\item[g)] $\displaystyle\int_{0}^{\theta}\mbox{\rm Sn}^p\varphi\;\mbox{\rm Cs}^{2n-1}\varphi\;d\varphi=\dfrac{\mbox{\rm Sn}^{p+1}\theta}{p+1}$;
		\item[h)] $\displaystyle\int_{0}^{T}\mbox{\rm Sn}^p\varphi\;\mbox{\rm Cs}^q\varphi\;d\varphi=0$, if either $p$ or $q$ are odd;
		\item[i)] $\displaystyle\int_{0}^{T}\mbox{\rm Sn}^p\varphi\;\mbox{\rm Cs}^q\varphi\;d\varphi=\frac{2}{\sqrt{n^{p+1}}}\dfrac{\Gamma\left(\frac{p+1}{2}\right)\Gamma\left(\frac{q+1}{2n}\right)}{\Gamma\left(\frac{p+1}{2}+\frac{q+1}{2n}\right)}$, if both $p$ and $q$ are even;
	\end{itemize}
\end{proposition}

Applying the generalized polar coordinates in system \eqref{FNfamilia} yields:

\begin{equation*}%\label{FNeqPolaresGen}
	\begin{array}{lcr}
		\dot{r}=r^{n}\left(\mu\mbox{Cs}^{n-1}\theta-(n-1)\mbox{Cs}^{2n-1}\theta\,\mbox{Sn}\,\theta+P(r,\theta)\right),\\
		\dot{\theta}=-nr^{n-1}\left(1-(n-1)\mbox{Sn}^2\theta+Q(r,\theta)\right),
	\end{array}
\end{equation*}
where $P(r,\theta)$ and $Q(r,\theta)$ are given by:
\begin{eqnarray}
P(r,\theta)&=&\displaystyle\sum_{i+nj\geqslant n+1}\tilde{a}_{ij}r^{i+nj-n}\mbox{Cs}^{i+2n-1}\theta\,\mbox{Sn}^j\theta\nonumber\\
&&+\displaystyle\sum_{i+nj\geqslant 2n}\tilde{b}_{ij}r^{i+nj-2n+1}\mbox{Cs}^{i}\theta\,\mbox{Sn}^{j+1}\theta,\nonumber\\
Q(r,\theta)&=&\displaystyle\sum_{i+nj\geqslant n+1}\tilde{a}_{ij}r^{i+nj-n}\mbox{Cs}^{i}\theta\,\mbox{Sn}^{j+1}\theta\nonumber\\
&&-\displaystyle\sum_{i+nj\geqslant 2n}\frac{\tilde{b}_{ij}}{n}r^{i+nj-2n+1}\mbox{Cs}^{i+1}\theta\,\mbox{Sn}^{j}\theta.\nonumber
\end{eqnarray}
We then obtain the following differential equation:
\begin{equation}\label{FNeqRgen}
	\dfrac{dr}{d\theta}=R(r,\theta)=-r\dfrac{\left(\mu\mbox{Cs}^{n-1}\theta-(n-1)\mbox{Cs}^{2n-1}\theta\,\mbox{Sn}\,\theta+P(r,\theta)\right)}{n\left(1-(n-1)\mbox{Sn}^2\theta+Q(r,\theta)\right)}.
\end{equation}
As $1-(n-1)\mbox{Sn}^2\theta=\mbox{Cs}^{2n}\theta+\mbox{Sn}^2\theta$, the denominator of \eqref{FNeqRgen} does not vanish at $r=0$. Thus $R(r,\theta)$ is an analytic function in a neighborhood of $r=0$. Let $\tilde{r}(\theta,r_0)$ be the solution of \eqref{FNeqRgen} with initial condition $\tilde{r}(0,r_0)=r_0$. We can expand it as a power series, i.e.
$$\tilde{r}(\theta,r_0)=\displaystyle\sum_{k=1}^{\infty}v_k(\theta)r_0^k,$$
and define the so-called \emph{focal values}. 

\begin{definition}
	Consider function $\tilde{r}(\theta,r_0)$ associated to equation \eqref{FNeqRgen} defined above. The \emph{displacement function} is defined by $$d(r_0)=\tilde{r}(T,r_0)-r_0.$$
	The \emph{focal values} are the coefficients of the power series expansion of $d(r_0)$ given by $v_1(T)-1$ and $v_k(T)$ for $k\geqslant 2$.
\end{definition}

The origin is a center for system \eqref{FNfamilia} if and only if all the focal values are null. In the literature, for the non-degenerate singular point, it is well known that the first non-zero focal value is the coefficient of an odd power of $r_0$ (see \cite[Lemma 5, page 243]{AndronovB} or \cite[Proposition 3.1.4, page 94]{Romanovski}). In \cite{LiuLi3Ord}, the authors studied system \eqref{FNfamilia} for $n=2$ and concluded that the first non-zero focal value is the coefficient of an even power of $r_0$. As a matter of fact, a more general result is true.

\begin{proposition}\label{PropoFNPrimeironnulogen}
	For system \eqref{FNfamilia} with monodromic nilpotent singular points having Andreev number $n$, the first non-zero focal value is the coefficient of a power of $r_0$ whose parity is the same as the parity of $n$.
\end{proposition}

This proposition is proved in \cite[pages 50 through 54]{Lyapunov} and will be useful in the next sections. Now we provide some results that will lead to the proof of Theorem \ref{TeoFocoNilpotenteForte}. 

\begin{lemma}
	The coefficient $v_1$ in function $\tilde{r}(\theta,r_0)$ associated to equation \eqref{FNeqRgen} is given by:
	\begin{equation}\label{FNeqv1gen}
		\exp\left(-\displaystyle\int_{0}^{\theta}\dfrac{\left(\mu\mbox{Cs}^{n-1}\varphi-(n-1)\mbox{Cs}^{2n-1}\varphi\,\mbox{Sn}\,\varphi\right)}{n\left(1-(n-1)\mbox{Sn}^2\varphi\right)}d\varphi\right).
	\end{equation}
\end{lemma}
\noindent\textbf{Proof: }Substituting the power series expansion of $\tilde{r}(\theta,r_0)$ in the equation \eqref{FNeqRgen} and comparing the coefficients of $r_0$ yields the result.\qed

\begin{lemma}\label{Lemav1gen}
	Consider function $\tilde{r}(\theta,r_0)$ associated to equation \eqref{FNeqRgen}. Then $v_1(T)=1$ for $n$ even and $v_1(T)=e^{-\mu A^2}$ for $n$ odd, where $A\neq 0$.
\end{lemma}
\noindent\textbf{Proof: }Using the properties described in Proposition \ref{propoPropriedadesTrigo}, we compute the expression \eqref{FNeqv1gen} for $\theta=T$:
$$v_1(T)=\exp\left(-\mu\displaystyle\int_{0}^{T}\dfrac{\left(\mbox{Cs}^{n-1}\theta\right)}{n\left(1-(n-1)\mbox{Sn}^2\theta\right)}d\theta\right).$$

Let $G(\theta)$ be the integrand of the above expression. For $n$ odd, $G(\theta)$ is a non-negative even function. Thus, the above integral is not zero, and then $v_1(T)=e^{-\mu A^2}$ for some $A\neq 0$. If $n$ is even, $G(\theta)$ is an odd $T$-periodic function. Therefore $v_1(T)=1$.\qed
\medskip

\noindent\textbf{Proof of Theorem \ref{TeoFocoNilpotenteForte}: }Consider system \eqref{FNeq1} under the hypothesis of Theorem \ref{TeoFocoNilpotenteForte}. By Proposition \ref{PropoPreformMono}, it can be written as system \eqref{FNfamilia}. By Remark \ref{obsbeta=n-1}, we have $\mu\neq0$ in canonical form \eqref{FNfamilia}. Lemma \ref{Lemav1gen} implies that $v_1(T)\neq 1$ for $\mu\neq 0$. Therefore, the origin cannot be a center, and so, it must be a nilpotent focus. \qed

\begin{remark}
It is possible to obtain similar properties and results if instead of the generalized polar coordinates, one make the change of variables  $x=r\cos\theta,\;y=r^n\sin\theta$, which is exactly the approach described by the authors of \cite{LiuLiNew1}. The choice of one change of variables over the other has implications in the simplicity of computations made by softwares of symbolic mathematics, such as Maple or Mathematica for instance. We also remark that the same results are obtained for family	\eqref{FNfamilia2}, aside the exact expression for $v_1(T)$.
\end{remark}

\section{Formal Series and Nilpotent singular points}

Let $H^{(n)}_{(2,1)}$ denote the vector space of homogeneous polynomials of degree $n$ in two variables. A base for this vector space is the collection of monomials $x^iy^j$ where $i,j\geqslant 0$ and $i+j=n$. For a given polynomial $p\in H^{(n)}_{(2,1)}$, we denote by $\langle p\rangle$ the vector subspace spanned by $p$. Consider the following linear transformation:

$$\begingroup
\setlength\arraycolsep{0pt}
T_n\colon\begin{array}[t]{c >{{}}c<{{}} c}
	H^{(n)}_{(2,1)} & \to & H^{(n)}_{(2,1)} \\
	\noalign{\medskip} p & \mapsto & y\dfrac{\partial p}{\partial x} .
\end{array}
\endgroup$$

\begin{lemma}\label{LemaFNKerTn}
	The kernel of linear transformation $T_n$ is given by $\ker T_n=\langle y^n\rangle$.
\end{lemma}
\noindent\textbf{Proof:} Clearly $T_n(y^n)=0$. Let $p=\sum_{j+k=n}p_{j,k}x^jy^k$. If $p\in\ker T_n$, then $T_n(p)=0$, which equates to:
$$\sum_{j+k=n}(j+1)p_{j+1,k-1}x^{j}y^{k}=0.$$
The coefficient of monomial $x^jy^k$ in the above expression is given by $(j+1)p_{j+1,k-1}$, which implies that $p_{j,k}=0$ for $j>0$. Thus $p=p_{0,n}y^n\in\langle y^n\rangle$.
\qed

\begin{lemma}\label{LemaFNTp+q}
	For every $q\in H^{(n)}_{(2,1)}$, there is a choice of $p\in H^{(n)}_{(2,1)}$ such that $T_n(p)+q\in\langle x^n\rangle$.
\end{lemma}
\noindent\textbf{Proof: }It is sufficient to prove $H^{(n)}_{(2,1)}=\langle x^n\rangle\oplus \text{Im} T_n$. By Lemma \ref{LemaFNKerTn}, the codimension of Im$T_n$ is 1. Since $\langle x^n\rangle$ has dimension 1, we only need to prove that the intersection of both vector subspaces is $\{0\}$. We have $T_n(p)\in\langle x^n\rangle$ if and only if $y\dfrac{\partial p}{\partial x}=\alpha x^n$, which is only possible when $\dfrac{\partial p}{\partial x}=\alpha=0$. Therefore $T_n(p)\in\langle x^n\rangle$ if and only if $p\in\langle y^n\rangle=\ker T_n$. The result follows. \qed

\bigskip

Given a function $F:\mathbb{R}^2\to\mathbb{R}$, we denote by $F_n$ the homogeneous part of degree $n$ in its power series expansion.

\begin{proposition}\label{PropoFNXH=xn}
	Consider a vector field $X$ associated to system \eqref{FNeq1} having a nilpotent singular point. Then there exists a formal series $H(x,y)=y^2+\sum_{n=3}^{\infty}H_n(x,y)$ such that $XH=\sum_{n\geqslant 4}\omega_nx^n$.
\end{proposition}
\noindent\textbf{Proof: }The Lie derivative $XH$ is given by:
\begin{small}
	\begin{eqnarray}
		XH&=&\left(y+P_2+P_3+\dots\right)\left(\dfrac{\partial H_3}{\partial x}+\dfrac{\partial H_4}{\partial x}+\dots\right)\nonumber\\
		&&+\left(Q_2+Q_3+\dots\right)\left(2y+\dfrac{\partial H_3}{\partial y}+\dfrac{\partial H_4}{\partial y}+\dots\right).\nonumber
	\end{eqnarray}
\end{small}

Note that $y\dfrac{\partial H_n}{\partial x}$ is a homogeneous polynomial of degree $n$. Hence, rewriting the above expression by organizing the homogeneous terms of degree $n$ in brackets, we have that
\begin{small}
	\begin{eqnarray}
		XH &=&\left(y\dfrac{\partial H_3}{\partial x}+2yQ_2\right)+\left(y\dfrac{\partial H_4}{\partial x}+F_4\right)+\left(y\dfrac{\partial H_5}{\partial x}+F_5\right)+\dots\nonumber\\
		&=&\sum_{n\geqslant 3}\left(y\dfrac{\partial H_n}{\partial x}+F_n\right)=\sum_{n\geqslant 3}T_n(H_n)+F_n,\nonumber
	\end{eqnarray}
\end{small}
where $F_n\in H^{(n)}_{(2,1)}$ are obtained by the homogeneous parts of $P,Q$ and $H$ of degree less than $n$. By Lemma \ref{LemaFNTp+q}, we can choose $H_n$ such that $T_n(H_n)+F_n=\omega_nx^n$, $\omega_n\in\mathbb{R}$. Moreover, for $n=3$, $F_3=2yQ_2\notin\langle x^3\rangle$ which means that there is a choice of $H_3$ such that $T_3(H_3)=-2yQ_2$. Then, by making the suitable choices, $XH=\sum_{n\geqslant 4}\omega_nx^n$.\qed

\bigskip

The following result is proven by following essentially the same steps in the last proof.

\begin{proposition}\label{PropoFNFatorInverso}
	Consider a vector field $X$ associated to system \eqref{FNeq1} having a nilpotent singular point. Then there exists a formal series $V(x,y)$ such that $XV-V\mbox{\rm div}X=\sum_{n\geqslant 1}\Lambda_nx^n$.
\end{proposition}

\noindent\textbf{Proof: }Consider $V(x,y)=\sum_{n=0}^{\infty}V_n(x,y)$. Then
\begin{small}
	\begin{eqnarray}
		XV-V\mbox{div}X &=&\left(y+P_2+P_3+\dots\right)\left(\dfrac{\partial V_1}{\partial x}+\dfrac{\partial V_2}{\partial x}+\dots\right)\nonumber\\
		&&+\left(Q_2+Q_3+\dots\right)\left(\dfrac{\partial V_1}{\partial y}+\dfrac{\partial V_2}{\partial y}+\dots\right)\nonumber\\
		&&-(V_0+V_1+V_2+\dots)\left(\dfrac{\partial P}{\partial x}+\dfrac{\partial Q}{\partial y}\right).\nonumber
	\end{eqnarray}
\end{small}
Rewriting the above expression by organizing the homogeneous terms, we have:
\begin{small}
	\begin{eqnarray}
		XV-V\mbox{div}X =\sum_{n\geqslant 1}\left(y\dfrac{\partial V_n}{\partial x}+F_n\right)=\sum_{n\geqslant 1}T_n(V_n)+F_n,\nonumber
	\end{eqnarray}
\end{small}
where $F_n\in H^{(n)}_{(2,1)}$ are obtained by the homogeneous parts of $P,Q$ and $V$ of degree less than $n$. By the same arguments as in the proof of Proposition \ref{PropoFNXH=xn}, the result follows.\qed

\bigskip

It is useful to notice that formal series $H(x,y)$ in Proposition \ref{PropoFNXH=xn} is constructed in an attempt to satisfy formally $XH\equiv 0$, that is, to be a \emph{formal first integral} for vector field $X$. The quantities $\omega_k$ which are non-zero pose an obstacle to this attempt.

Proposition \ref{PropoFNXH=xn} is analogous to a well known result for non-degenerate singular point, except that in this case, $j^2H(0)=x^2+y^2$ and $XH$ has only even powers of $(x^2+y^2)$, (see, for instance, \cite{Llibre,Romanovski}). Moreover, any obstruction to $H$ being a formal first integral implies the existence of a focus in the singular point, due to the fact that $j^nH(0)$, for $n$ big enough, is a \emph{Lyapunov function} \cite{Lyapunov}. In the nilpotent case, this is not necessarily true. The expressions of $H(x,y)$, in particular $j^2H(0)=y^2$, and $XH$ given in Proposition \ref{PropoFNXH=xn} do not immediately imply that $j^nH(0)$, for $n$ big enough, will be a Lyapunov function. However, the following result provides a useful tool to study the Center problem when the first obstruction to $XH=0$ is an even power of $x$.

\begin{proposition}\label{PropoFNXHNFoco}
	Let $X$ be the vector field associated to system \eqref{FNeq1} having a nilpotent singular point and $H$ be a formal series as in Proposition \ref{PropoFNXH=xn}. If there exists $n\in\mathbb{N}$ such that $j^{2n}XH(0)=\omega_{2n}x^{2n}$ with $\omega_{2n}\neq 0$, then the origin cannot be a center.
\end{proposition}
\noindent\textbf{Proof: }Let $H(x,y)=y^2+\sum_{j=3}^{\infty}H_j(x,y)$ be the formal series given by Proposition \ref{PropoFNXH=xn}. Consider its $2n$-jet, i.e.: $\tilde{H}(x,y)=y^2+\sum_{j=3}^{2n}H_j(x,y)$. Thus, $\tilde{H}$ is an analytic function such that $j^{2n}X\tilde{H}(0)=\omega_{2n}x^{2n}$.

Without loss of generality suppose $\omega_{2n}<0$. The same conclusions will follow for $\omega_{2n}>0$ considering vector field $-X$.

Since $j^2\tilde{H}(0)=y^2$ and $j^{2n}X\tilde{H}(0)=\omega_{2n}x^{2n}$, there is a neighborhood $U$ of the origin such that $\tilde{H}\geqslant 0$ and $X\tilde{H}\leqslant 0$. Let $p\in U$ and $\gamma_p(t)$ the trajectory of system \eqref{FNeq1} with initial point $\gamma_p(0)=p$. Note that $\tilde{H}$ is non-increasing along $\gamma_p$. If $\gamma_p$ is periodic and wholly contained in $U$, there is $T\in\mathbb{R}$ such that $\gamma_p(0)=\gamma_p(T)$ and thus $\tilde{H}\circ\gamma_p(0)=\tilde{H}\circ\gamma_p(T)$. By continuity of $\tilde{H}$, we have $X\tilde{H}\vert_{\gamma_p}\equiv 0$. Therefore, any closed orbits of $X$ inside $U$ are contained in the level sets of $\tilde{H}$. Thus, if the origin is a center, there is a neighborhood of it for which $X\tilde{H}$ is null. But this contradicts the hypothesis.\qed

\begin{remark}
	Using Proposition \ref{PropoFNXHNFoco} we can give another proof of Theorem \ref{TeoFocoNilpotenteForte}. For system \eqref{FNfamilia} with $n$ odd, we compute the first non-zero $\omega_k$ in $XH$, which is $\omega_{3n-1}=4n\mu(\mu^2+1)$ (the computations are omitted for the sake of simplicity). Since $3n-1$ is even, we cannot have a center at the origin, and by monodromy, the singular point is a nilpotent focus.
\end{remark}

\section{Shortcomings of the inverse integrating factor method for systems \eqref{FNeq1nhomo2n-1}}\label{secImpossibilidade}

\noindent\textbf{Proof of Theorem \ref{TeoFNimpossibilidadeGarcia}: }Consider $(1,n)$-quasi-\-homo\-geneous systems \eqref{FNeq1nhomo2n-1} having a monodromic nilpotent singular point with Andreev number $n$ odd at the origin and maximum weighted degree $2n-1$. Let  $V(x,y)=\sum_{i+j\geqslant 0}q_{ij}x^iy^j$ be a formal inverse integrating factor for system \eqref{FNeq1nhomo2n-1}. Note that if $q_{00}\neq 0$, i.e. $V(0,0)\neq 0$, then there exists a formal first integral $H(x,y)$ for system \eqref{FNeq1nhomo2n-1} (see \cite{GarciaFiI}) and $j^2H(0)=y^2$ (see \cite{Chavarriga}). Hence, Proposition \ref{PropoFNXH=xn} allow us to implement an algorithm to search for a formal first integral $H(x,y)$ to system \eqref{FNeq1nhomo2n-1}. More precisely, we will see in Subsection \ref{Computation} that the obstructions for the existence of such formal first integral are the quantities $\omega_{3n},\dots,\omega_{2(2n-1)}$ given by:

$$\omega_{3n}=2a_{n+1,0},$$
\begin{equation}\label{FNeqVkImpossibilidade}
	\omega_{3n+k}\equiv 2a_{n+1+k,0}\mod\langle \omega_{3n},\dots,\omega_{3n+k-1}\rangle,\;k=1,\dots,n-2.
\end{equation}

Therefore the conditions on the parameters for formal integrability are $a_{n+1,0}=\dots=a_{2n-1,0}=0$. Under these conditions, system \eqref{FNeq1nhomo2n-1} is invariant by transformation $x\to x,\;y\to-y,\;t\to-t$, and so the origin is a nilpotent center.

Now we use Proposition \ref{PropoFNFatorInverso} to find conditions for the system to have a formal inverse integrating factor with $q_{00}=0$. The obstructions, computed in Subsection \ref{Computation}, are now the quantities $\Lambda_{3n},\dots,\Lambda_{2(2n-1)}$ given by:

$$\Lambda_{3n}=-\dfrac{a_{n+1,0}q_{02}}{n},$$
\begin{equation}\label{FNeqLambdakImpossibilidade}
	\Lambda_{3n+k}\equiv -\dfrac{(k+1)a_{n+1+k,0}q_{02}}{n}\mod\langle \Lambda_{3n},\dots,\Lambda_{2(2n-1)}\rangle,
\end{equation}
for $k=1,\dots,n-2$. Moreover, the coefficients $q_{02},q_{n1},q_{2n,0}$ satisfy the equations $q_{2n,0}=\frac{q_{02}}{n}$ and $q_{n1}=0$. For $q_{02}\neq 0$, we obtain the same center conditions of the case $q_{00}\neq 0$ and so, no new nilpotent centers are detected. For $q_{02}=0$, the $(1,n)$-quasi-homogeneous term of degree $2n$ in $V(x,y)$ is null and Theorem \ref{TeoGarciaFiI2} cannot be applied to solve the center problem. The conclusion is Theorem \ref{TeoFNimpossibilidadeGarcia}. \qed

\begin{remark}
	Applying Theorem \ref{TeoGarciaFiI2} to systems \eqref{FNeq1nhomo2n-1} with $n$ even also yields the center conditions $a_{k0}=0$ for $k=n+1,\dots,2n-1$.
\end{remark}

Since the inverse integrating factor method cannot detect all center conditions for system \eqref{FNeq1nhomo2n-1} with $n$ odd, it is necessary to search for other approaches in order to solve the center problem. System \eqref{FNeq1nhomo2n-1} is in the class of $(1,n)$-quasi-homogeneous systems studied in \cite{GasullTorre}. In this work, the authors propose a way of computing the focal values in a recurrent way using the $(1,n)$-quasi-homogeneous terms of the system after performing the generalized polar coordinate change. We will use this method to prove Proposition \ref{PropoFNcondcentro} which presents a necessary condition for system \eqref{FNeq1nhomo2n-1} to have a center that was not obtainable by the inverse integrating factor method.

\begin{lemma}\label{LemmaV3}
For system \eqref{FNfamilia2}, with $\tilde{\mu}=0$ and $n$ odd, the focal value $v_3(T)$ is given by:
\begin{eqnarray}\label{FNeqv3}
v_3(T)&=&\epsilon_3\left((n+2)(\hat{a}_{n+2,0}+\hat{a}_{n+1,0}(\hat{a}_{11}+2\hat{b}_{02}))\right.\nonumber\\
&&\left.+\hat{b}_{n+1,1}+\hat{b}_{n,1}(\hat{a}_{11}+2\hat{b}_{02})\right),
\end{eqnarray}
where $\epsilon_3$ is a non-zero constant.
\end{lemma}
\noindent\textbf{Proof: }The generalized polar coordinates transform \eqref{FNfamilia2}, with $\tilde{\mu}=0$, into the following system
\begin{equation}\label{FNeqPolar}
	\begin{array}{lcr}
		\dot{r}=r^{n+1}\sum_{k\geqslant 0}r^kR_{3n+k}(\theta),\\
		\dot{\theta}=r^{n-1}\left(1+\sum_{k\geqslant 1}r^k\Theta_{2n+k}(\theta)\right),\\
	\end{array}
\end{equation}
where $R_{k}(\theta)$ and $\Theta_k(\theta)$ are $(1,n)$-quasi-homogeneous polynomials of weighted degree $k$ in the variables $\mbox{Cs}\,\theta$, $\mbox{Sn}\,\theta$. Using the method described in \cite[Theorem 1]{GasullTorre}, we have that
$$v_3(T)=\varepsilon_3\int_{0}^{T}h_1r^{2n+1}R_{3n}(\theta)+r^{2n+2}R_{3n+1}(\theta)\,d\theta,$$
for which $h_1=-r\int_{0}^{\theta}\Theta_{2n+1}'(\varphi)+(2n+1)R_{3n}(\varphi)\,d\varphi$ and $\varepsilon_3$ is a non-zero constant. Thus, we only need the expressions of $R_{3n}$, $R_{3n+1}$ and $\Theta_{2n+1}$ in order to compute $v_3(T)$. For system \eqref{FNeqPolar} those are given by:
\begin{eqnarray}
R_{3n}(\theta)&=&\hat{a}_{n+1,0}\mbox{Cs}^{3n}\theta+(\hat{a}_{11}+\hat{b}_{2n,0})\mbox{Cs}^{2n}\theta\,\mbox{Sn}\,\theta\nonumber\\
&&+\hat{b}_{n,1}\mbox{Cs}^{n}\theta\,\mbox{Sn}^2\theta+\hat{b}_{02}\mbox{Sn}^3\theta,\nonumber\\
R_{3n+1}(\theta)&=&\hat{a}_{n+2,0}\mbox{Cs}^{3n+1}\theta+(\hat{a}_{21}+\hat{b}_{2n+1,0})\mbox{Cs}^{2n+1}\theta\mbox{Sn}\,\theta\nonumber\\
&&+\hat{b}_{n+1,1}\mbox{Cs}^{n+1}\theta\mbox{Sn}^{2}\theta+\hat{b}_{12}\mbox{Cs}\,\theta\mbox{Sn}^3\theta,\nonumber\\
\Theta_{2n+1}(\theta)&=&\hat{b}_{2n,0}\mbox{Cs}^{2n+1}\theta+(\hat{b}_{n,1}-n\hat{a}_{n+1,0})\mbox{Cs}^{n+1}\theta\,\mbox{Sn}\,\theta\nonumber\\
&&+(\hat{b}_{02}-n\hat{a}_{11})\mbox{Cs}\,\theta\,\mbox{Sn}^2\theta.\nonumber
\end{eqnarray}
After extensive calculations, using Proposition \ref{propoPropriedadesTrigo}, the above expression simplifies to
\begin{eqnarray}
v_3(T)&=&\varepsilon_3\int_{0}^{T}(\hat{a}_{n+2,0}+\hat{a}_{n+1,0}(\hat{a}_{11}+2\hat{b}_{02}))\mbox{Cs}^{3n+1}\theta\,d\theta\nonumber\\
&&+\varepsilon_3\int_{0}^{T}
(\hat{b}_{n+1,1}+\hat{b}_{n,1}(\hat{a}_{11}+2\hat{b}_{02}))\mbox{Cs}^{n+1}\theta\,\mbox{Sn}^2\theta\,d\theta.\nonumber
\end{eqnarray}
It is easy to see that $\int_{0}^{T}\mbox{Cs}^{3n+1}\theta\,d\theta=(n+2)\int_{0}^{T}
\mbox{Cs}^{n+1}\theta\,\mbox{Sn}^2\theta\,d\theta$. Since $n$ is odd, $\int_{0}^{T}
\mbox{Cs}^{n+1}\theta\,\mbox{Sn}^2\theta\,d\theta$ is a positive number. Setting $$\epsilon_3:=\varepsilon_3\int_{0}^{T}
\mbox{Cs}^{n+1}\theta\,\mbox{Sn}^2\theta\,d\theta,$$ yields the expression \eqref{FNeqv3}. \qed

\bigskip

\noindent\textbf{Proof of Proposition \ref{PropoFNcondcentro}: }Consider system \eqref{FNeq1nhomo2n-1} with $n$ odd. By Lemma \ref{Lemav1gen} and Proposition \ref{PropoFNPrimeironnulogen}, we know that $v_1(T)=v_2(T)=0$. Using Lemma \ref{LemmaV3}, we have that $v_3(T)=\varepsilon_3(a_{11}a_{n+1,0}+a_{n+2,0})$. Thus a necessary condition for the origin of system \eqref{FNeq1nhomo2n-1} to be a nilpotent center is $a_{n+2,0}=-a_{11}a_{n+1,0}$. We have proven Proposition \ref{PropoFNcondcentro}. \qed

\subsubsection{Computation of the obstructions \eqref{FNeqVkImpossibilidade} and \eqref{FNeqLambdakImpossibilidade}}\label{Computation}

Denote by $X$ the vector field associated to system \eqref{FNeq1nhomo2n-1}. To obtain the equations \eqref{FNeqVkImpossibilidade}, it is sufficient to analyze the coefficients of the power series expansion of $XH$. We write $H(x,y)=y^2+\sum_{i+j\geqslant 3}p_{ij}x^iy^j$ and define by convention $p_{02}=1, p_{ij}=0$ if ${i+j\leqslant 2, j\neq 2}$ or if $i,j<0$. Then:
\begin{small}
	\begin{eqnarray}
		XH&=&-\sum_{i+j\geqslant 3}ip_{ij}x^{i-1}y^{j+1}+\sum_{i+j\geqslant 3}ip_{ij}x^{i-1}y^j\displaystyle\sum_{k=1}^{n-1}a_{k1}x^{k}y\nonumber\\
		&&+\sum_{i+j\geqslant 3}ip_{ij}x^{i-1}y^j\displaystyle\sum_{k=n+1}^{2n-1}a_{k0}x^{k}+2x^{2n-1}y+\sum_{i+j\geqslant 3}jp_{ij}x^{i+2n-1}y^{j-1}.\nonumber
	\end{eqnarray}
\end{small}

The coefficient $W_{kl}$ of monomial $x^ky^l$ in $XH$ is given by: 
\begin{eqnarray}\label{Wkl}
	W_{kl}&=&-(k+1)p_{k+1,l-1}+\sum_{i=k-n+2}^{k}ip_{i,l-1}a_{k-i+1,1}\nonumber\\
	&&+\sum_{i=k-2n+2}^{k-n}ip_{i,l}a_{k-i+1,0}+(l+1)p_{k-2n+1,l+1}.
\end{eqnarray}
Thus, by Proposition \ref{PropoFNXH=xn}, $W_{kl}=0$ for $l\neq 0$ and
$$\omega_k=W_{k0}=\sum_{i=k-2n+2}^{k-n}ip_{i,0}a_{k-i+1,0}+p_{k-2n+1,1}.$$
From \eqref{Wkl}, since $W_{k1}=0$, for $k\leqslant n$ we have:
$$(k+1)p_{k+1,0}=\sum_{i=2}^{k}ip_{i,0}a_{k-i+1,1}.$$
Hence, for $k=1$, we conclude that $p_{20}=0$. This implies that $p_{30}=0$ for $k=2$. Repeating this procedure, we obtain that $p_{k0}=0$ for $k=0,\dots n+1$.
%$$(n+1)p_{n+1,0}=\sum_{i=2}^{n}ip_{i,0}a_{n-i+1,1}$$
Now, from $W_{k2}=0$, for $k\leqslant n$ it follows that:
$$(k+1)p_{k+1,1}=\sum_{i=1}^{k}ip_{i,1}a_{k-i+1,1}.$$
By recalling that $p_{11}=0$ we have, as the previous case, that $p_{k,1}=0$ for $k=0,\dots,n+1$.
%$$(n+1)p_{n+1,1}=\sum_{i=2}^{n}ip_{i,1}a_{n-i+1,1}$$
Therefore:
$$\omega_{3n}=\sum_{i=n+2}^{2n}ip_{i,0}a_{3n-i+1,0}.$$
Since $W_{k1}=0$, for $n+1\leqslant k< 2n-1$, we have:
$$(k+1)p_{k+1,0}=\sum_{i=n}^{k}ip_{i,0}a_{k-i+1,1},$$
which imply that $p_{k0}=0$ for $k=n+2,\dots,2n-1$. Now for $k=2n-1$, by \eqref{Wkl}, $W_{2n-1,1}=0$ is equivalent to
$$(2n)p_{2n,0}=\sum_{i=n+1}^{2n-1}ip_{i,0}a_{2n-i,1}+2=2\Rightarrow p_{2n,0}=\dfrac{1}{n}.$$
Thus, we obtain 
$$\omega_{3n}=2a_{n+1,0},$$
and $\omega_{3n}=0$ if and only if $a_{n+1,0}=0$.
Computing the quantities $\omega_{3n+k}$ for $k=1,\dots,n-2$, we have:
\begin{eqnarray}\label{omega3n+k}
	\omega_{3n+k}&=&\sum_{i=2n}^{2n+k}ip_{i,0}a_{3n+k-i+1,0}+p_{n+k+1,1}=2a_{n+1+k,0}+p_{n+k+1,1}\nonumber\\
	&&+\sum_{i=2n+1}^{2n+k}ip_{i,0}a_{3n+k-i+1,0}.\nonumber
\end{eqnarray}
Suppose, by induction, that 
$$\omega_{3n+j}\equiv 2a_{n+1+j,0} \mod\langle \omega_{3n},\dots, \omega_{3n+j-1}\rangle,$$
and 
$$ \mbox{ and } p_{n+j,1}\equiv 0\mod\langle \omega_{3n},\dots, \omega_{3n+j-1}\rangle,$$ 
for all $j$ such that ${0\leqslant j\leqslant k-1}$. Thus, $\omega_{3n}=\dots=\omega_{3n+k-1}=0$ if and only if $a_{n+1}=\dots=a_{n+k}=0$. And so \eqref{omega3n+k} becomes:
$$\omega_{3n+k}\equiv 2a_{n+1+k,0}+p_{n+k+1,1}\mod \langle \omega_{3n},\dots, \omega_{3n+k-1}\rangle.$$
By \eqref{Wkl}, we obtain:
$$(n+k+1)p_{n+k+1,1}\equiv\sum_{i=n}^{n+k}ip_{i,1}a_{n+k-i+1,1}\mod \langle \omega_{3n},\dots, \omega_{3n+k-1}\rangle.$$
Therefore $p_{n+k+1,1}\equiv 0\mod \langle \omega_{3n},\dots \omega_{3n+k-1}\rangle$ and consequently
$$\omega_{3n+k}\equiv 2a_{n+1+k,0}\mod \langle \omega_{3n},\dots, \omega_{3n+k-1}\rangle.$$
Hence, \eqref{FNeqVkImpossibilidade} hold.

The equations \eqref{FNeqLambdakImpossibilidade} are obtained by performing a similar analysis in the following expression:
\begin{small}
	\begin{eqnarray}
		XV-V\mbox{div}X&=&-\sum_{i+j\geqslant 0}iq_{ij}x^{i-1}y^{j+1}+\sum_{i+j\geqslant 0}iq_{ij}x^{i-1}y^j\displaystyle\sum_{k=1}^{n-1}a_{k1}x^{k}y\nonumber\\
		&&+\sum_{i+j\geqslant 0}iq_{ij}x^{i-1}y^j\displaystyle\sum_{k=n+1}^{2n-1}a_{k0}x^{k}+\sum_{i+j\geqslant0}jq_{ij}x^{i+2n-1}y^{j-1}\nonumber\\
		&&-\sum_{i+j\geqslant 0}q_{ij}x^iy^j\displaystyle\sum_{k=1}^{n-1}ka_{k1}x^{k-1}y\nonumber\\
		&&-\sum_{i+j\geqslant 0}q_{ij}x^iy^j\displaystyle\sum_{k=n+1}^{2n-1}ka_{k0}x^{k-1}.\nonumber
	\end{eqnarray}
\end{small}
Again, we define $q_{ij}=0$ for $i+j<0$ or $i,j<0$. Analogous to the previous argument, the coefficient of monomial $x^ky^l$ above is given by:

\begin{eqnarray}\label{Lkl}
	L_{kl}=-(k+1)q_{k+1,l-1}+\sum_{i=k-n+2}^{k}(2i-k-1)q_{i,l-1}a_{k-i+1,1}+\nonumber\\
	+\sum_{i=k-2n+2}^{k-n}(2i-k-1)q_{il}a_{k-i+1,0}+(l+1)q_{k-2n+1,l+1}.
\end{eqnarray}

By Proposition \ref{PropoFNFatorInverso}, $L_{kl}=0$ for $l\neq 0$, and
\begin{equation}\label{Lambda3n}
	\Lambda_{k}=L_{k0}=\sum_{i=k-2n+2}^{k-n}(2i-k-1)q_{i0}a_{k-i+1,0}+q_{k-2n+1,1}.
\end{equation}
From \eqref{Lkl}, since $L_{k2}=0$, for $k\leqslant n-2$ we have:
\begin{eqnarray}
(k+1)q_{k+1,1}&=&\sum_{i=0}^{k}(2i-k-1)q_{i,1}a_{k-i+1,1}\nonumber\\
&=&-(k+1)q_{01}a_{k+1,1}+\sum_{i=1}^{k}(2i-k-1)q_{i,1}a_{k-i+1,1},\nonumber
\end{eqnarray}
which yields $q_{11}=-a_{11}q_{01}$. Now, for $k=1$, this implies that $q_{21}=-a_{21}q_{01}$. Continuing with this procedure we obtain $q_{k,1}=-a_{k,1}q_{01}$ for $k=1,\dots n-1$. To compute $q_{n1}$, we substitute $k=n-1$ and $l=2$ in \eqref{Lkl} and equate the expression to zero. That is:
\begin{equation}\label{qn1}
	nq_{n,1}=\sum_{i=1}^{n-1}(2i-n)q_{i,1}a_{n-i,1}=-q_{01}\sum_{i=1}^{n-1}(2i-n)a_{i,1}a_{n-i,1}=0.
\end{equation}

Now, substituting $k=n$ and $l=2$ in \eqref{Lkl}, from $L_{n2}=0$, we obtain:
\begin{eqnarray}
(n+1)q_{n+1,1}&=&\sum_{i=2}^{n-1}(2i-n-1)q_{i,1}a_{n-i+1,1}+(-n-1)q_{02}a_{n+1,0}\nonumber\\
&=&-(n+1)a_{n+1,0}q_{02},\nonumber
\end{eqnarray}
and therefore $q_{n+1,1}=-a_{n+1,0}q_{02}$. We now can compute some coefficients $q_{i0}$ using \eqref{Lkl}. For $k\leqslant n-2$ we have:
\begin{eqnarray}
(k+1)q_{k+1,0}&=&\sum_{i=0}^{k}(2i-k-1)q_{i,0}a_{k-i+1,1}\nonumber\\
&=&-(k+1)q_{00}a_{k+1,1}+\sum_{i=1}^{k}(2i-k-1)q_{i,0}a_{k-i+1,1}.\nonumber
\end{eqnarray}
Then we conclude that $q_{10}=-a_{11}q_{00}$. Analogously to the previous cases, this implies that $q_{k,0}=-a_{k,1}q_{00}$ for $k=1,\dots,n-1$. Substituting $k=n-1,n$ and $l=1$ in \eqref{Lkl} and equating the resulting expression to zero yields
$$nq_{n,0}=\sum_{i=1}^{n-1}(2i-n)q_{i,0}a_{n-i,1}=-q_{00}\sum_{i=1}^{n-1}(2i-n)a_{i,1}a_{n-i,1}=0,$$
and
$$(n+1)q_{n+1,0}=\sum_{i=2}^{n-1}(2i-n-1)q_{i,0}a_{n-i+1,1}-(n+1)q_{01}a_{n+1,0},$$
respectively. Thus $q_{n,0}=0$ and $q_{n+1,0}=-q_{01}a_{n+1,0}$. For $n+1\leqslant k< 2n-1$, by \eqref{Lkl} and the fact that $L_{k1}=0$, we have
\begin{eqnarray}
(k+1)q_{k+1,0}&=&\sum_{i=k-n+2}^{k}(2i-k-1)q_{i,0}a_{k-i+1,1}\nonumber\\
&&-(k+1)q_{01}a_{k+1,0}-q_{01}\sum_{i=1}^{k-n}(2i-k-1)a_{i1}a_{k-i+1,0}.\nonumber
\end{eqnarray}

Since, for $k=n+1$,
\begin{eqnarray}
(n+2)q_{n+2,0}&=&-q_{00}\sum_{i=3}^{n-1}(2i-n-2)a_{i1}a_{n-i+2,1}-(n+2)q_{01}a_{n+2,0}\nonumber\\
&=&-(n+2)q_{01}a_{n+2,0},\nonumber
\end{eqnarray}
by induction we have $q_{k0}=-a_{k0}q_{01}$ for $k=n+2,\dots,2n-1$. Then, substituting $k=2n-1$ and $l=1$ in \eqref{Lkl} and equating the resulting expression to zero, it follows that
\begin{eqnarray}
(2n)q_{2n,0}&=&-q_{01}\sum_{i=n+1}^{2n-1}(2i-2n)a_{i0}a_{2n-i,1}\nonumber\\
&&-q_{01}\sum_{i=1}^{n-1}(2i-2n)a_{i1}a_{2n-i,0}+2q_{02}\nonumber\\
&=&2q_{02},\nonumber
\end{eqnarray}
and therefore $q_{2n,0}=\dfrac{q_{02}}{n}$. By \eqref{qn1}, we have that $q_{n1}=0$ and the $(1,n)$-quasi-homogeneous term of weighted degree $2n$ of $V(x,y)$ is null if and only if $q_{02}=0$. For now on we assume $q_{02}\neq 0$.

Now we are able to compute $\Lambda_{3n}$ using \eqref{Lambda3n}. We have
\begin{eqnarray}
\Lambda_{3n}&=&(n-1)q_{2n,0}a_{n+1,0}-a_{n+1,0}q_{02}\nonumber\\
&&-q_{01}\sum_{i=n+2}^{2n-1}(2i-3n-1)a_{i0}a_{3n-i+1,0}\nonumber\\
&=&-\dfrac{a_{n+1,0}q_{02}}{n}.\nonumber
\end{eqnarray}

We now compute the quantities $\Lambda_{3n+k}$ for $1\leqslant k\leqslant n-2$:
\begin{equation}\label{Lambda3n+k}
	\Lambda_{3n+k}=\sum_{i=n+k+2}^{2n+k}(2i-3n+k-1)q_{i0}a_{3n+k-i+1,0}+q_{n+k+1,1}.
\end{equation}
Suppose, by induction, that 
$$\Lambda_{3n+j}\equiv -\dfrac{(j+1)a_{n+1+j,0}q_{02}}{n}\mod\langle \Lambda_{3n},\dots,\Lambda_{3n+j-1}\rangle,$$
and
$$q_{n+j+1,1}=-a_{n+j+1,0}q_{02}\mod\langle \Lambda_{3n},\dots,\Lambda_{3n+j-1}\rangle,$$
for every $j$ such that $0\leqslant j\leqslant k-1$. Since we are assuming $q_{02}\neq 0$, $\Lambda_{3n}=\dots=\Lambda_{3n+k-1}=0$ if and only if $a_{n+1,0}=\dots=a_{n+k,0}=0$. Then, equation \eqref{Lambda3n+k} becomes:
\begin{eqnarray}\label{Lambda3n+k2}
	\Lambda_{3n+k}&\equiv&\sum_{i=n+k+2}^{2n}(2i-3n-k-1)q_{i0}a_{3n+k-i+1,0}+q_{n+k+1,1}\nonumber\\
	&\equiv&(n-k-1)q_{2n0}a_{n+k+1,0}+q_{n+k+1,1}\nonumber\\
	&&-q_{01}\sum_{i=n+k+2}^{2n-1}(2i-3n-k-1)a_{i0}a_{3n+k-i+1,0}\nonumber\\
	&\equiv&\dfrac{(n-k-1)q_{02}a_{n+k+1,0}}{n}+q_{n+k+1,1}\nonumber\\
	&&\mod\langle \Lambda_{3n},\dots,\Lambda_{3n+k-1}\rangle.
\end{eqnarray}
Using \eqref{Lkl} again, we obtain:
\begin{eqnarray}
	(n+k+1)q_{n+k+1,1}&\equiv&\sum_{i=k+2}^{n-1}(2i-n-k-1)q_{i,1}a_{n+k-i+1,1}\nonumber\\
	& &+\sum_{i=n+1}^{n+k}(2i-n-k-1)q_{i,1}a_{n+k-i+1,1}\nonumber\\
	&&-(n+k+1)q_{02}a_{n+k+1,0}\nonumber\\
	&\equiv& -q_{01}\sum_{i=k+2}^{n-1}(2i-n-k-1)a_{i,1}a_{n+k-i+1,1}\nonumber\\
	& &-q_{02}\sum_{i=n+1}^{n+k}(2i-n-k-1)a_{i,0}a_{n+k-i+1,1}\nonumber\\
	&&-(n+k+1)q_{02}a_{n+k+1,0}\nonumber\\
	&\equiv& -(n+k+1)q_{02}a_{n+k+1,0}\nonumber\\
	&&\mod\langle \Lambda_{3n},\dots,\Lambda_{3n+k-1}\rangle.\nonumber
\end{eqnarray}
Thus
$$q_{n+k+1,1}=-a_{n+k+1,0}q_{02}\mod\langle \Lambda_{3n},\dots,\Lambda_{3n+k-1}\rangle,$$
and, by \eqref{Lambda3n+k2}, $\Lambda_{3n+k}\equiv -\dfrac{(k+1)a_{n+1+k,0}q_{02}}{n}\mod\langle \Lambda_{3n},\dots,\Lambda_{3n+k-1}\rangle$, and we have proven \eqref{FNeqLambdakImpossibilidade}.\qed

\section{Nilpotent centers for a system having a singular point with Andreev number $n=3$}\label{sec135}

\noindent\textbf{Proof of Theorem \ref{TeoFNCentro135}: }Considering system \eqref{FNeq13homo5}, by Theorem \ref{TeoFocoNilpotenteForte} the first necessary center condition for this system is $\mu=0$. Firstly, we search for integrable centers. We use Proposition \ref{PropoFNXH=xn} to construct a formal series $H(x,y)$, with $j^2H(0)=y^2$, such that $XH=\sum_{k\geqslant3}\omega_kx^k$. In the proof of Theorem \ref{TeoFNimpossibilidadeGarcia}, we have a general formula for some of the $\omega_k$. For the sake of completeness, we write them explicitly for system \eqref{FNeq13homo5}:
$$\omega_3=\dots=\omega_8=0,\;\omega_9=2a_{40},\;\omega_{10}=2(a_{40}a_{11}+a_{50}).$$
$$\omega_{11}=0,\;\omega_{12}\equiv\dfrac{2a_{11}a_{40}a_{21}}{7}\mod\langle \omega_{10}\rangle$$
The formal integrability implies $a_{40}=a_{50}=0$. Under those, system \eqref{FNeq13homo5} becomes
\begin{equation*}
	\begin{array}{lcr}
		\dot{x}=-y+a_{11}xy+a_{21}x^2y,\\
		\dot{y}=x^5,\\
	\end{array}
\end{equation*}
which is invariant by transformation ${x\to x},\;{y\to -y},$ ${t\to -t}$, and therefore the origin is a nilpotent center. We remark that $a_{40}=0\neq a_{50}$ is not a center condition by Proposition \ref{PropoFNXHNFoco}. Since not every nilpotent center is integrable, we need to continue our study.

We then turn to compute focal values for system \eqref{FNeq13homo5} with $\mu=0$. By Lemma \ref{Lemav1gen} and Proposition \ref{PropoFNPrimeironnulogen}, we get $v_1(T)=v_2(T)=0$. By \eqref{v3}, we have:

$$v_3(T)=\varepsilon_3(a_{11}a_{40}+a_{50})\left(\int_{0}^{T}\mbox{Cs}^{10}\theta\;d\theta\right),$$
where $\varepsilon_3$ is a non-zero constant. New necessary center condition is $a_{50}=-a_{11}a_{40}$, which cannot be obtained by the inverse integrating factor method. Assuming these conditions reduce our problem to $3$ parameters. We then analyze system \eqref{FNeq13homo5} as a limit of systems having non-degenerated centers using the algorithm described in paper \cite{LlibreLimite} and its previous versions. It consists in working with the following perturbation of system \eqref{FNeq13homo5} under $\mu=a_{11}a_{40}+a_{50}=0$:
\begin{equation*}
	\begin{array}{lcr}
		\dot{x}=-y+a_{11}xy+a_{21}x^2y+a_{40}x^4-a_{11}a_{40}x^5+\varepsilon P(x,y),\\
		\dot{y}=\varepsilon x+x^5+\varepsilon Q(x,y),\\
	\end{array}
\end{equation*}
where $j^1P(0)=j^1Q(0)=0$, and $\varepsilon>0$. We compute Lyapunov quantities for the above system, since it is non-degenerated for every $\varepsilon>0$, searching conditions on the parameters $a_{11},a_{21},a_{40}$ for which the origin of the perturbed system is a center for every positive $\varepsilon$ near $\varepsilon=0$. Denoting by $g_k$ the $k$-th Lyapunov constant, we have:
$$g_2=-\dfrac{8}{3}a_{40}(a_{11}+\tilde{q}_{{20}})\varepsilon+O(\varepsilon^2),$$
$$g_4=\frac{8}{189}a_{40}\left({53\,a_{11}}+{80\,\tilde{q}_{{20}}}\right)+O(\varepsilon),$$
where $\tilde{q}_{20}$ is the coefficient of $x^2$ in $Q(x,y)$. Assuming $a_{40}\neq 0$ yields $\tilde{q}_{20}=-a_{11}$ and consequentially ${g_4=\frac{8}{7}a_{40}a_{11}+O(\varepsilon)}$. Thus, for the origin to be a center of the perturbed system, we need $a_{11}=0$. Therefore $a_{50}=-a_{11}a_{40}=0$ is a necessary center condition for system \eqref{FNeq13homo5}. Those are also sufficient since under $a_{11}=a_{50}=\mu=0$ we have the following system:
\begin{equation*}
	\begin{array}{lcr}
		\dot{x}=-y+a_{21}x^2y+a_{40}x^4,\\
		\dot{y}=x^5,\\
	\end{array}
\end{equation*}
which is invariant by $y\to y,\;x\to -x,\;t\to -t$ and thus the origin is a nilpotent center. This proves Theorem \ref{TeoFNCentro135}. \qed

\begin{remark}
For $(1,n)$-quasi-homogeneous systems \eqref{FNeq1nhomo2n-1} of maximum weighted degree $2n-1$, if $a_{k1}=a_{k0}=0$ for all $k$ odd then the origin is a nilpotent center and the system is reversible. Moreover, this condition satisfies Proposition \ref{PropoFNcondcentro}.
\end{remark}

\medskip

\section{Acknowledgments}

The first author is partially supported by S\~ao Paulo Research Foundation (FAPESP) grants 19/10269-3 and 18/19726-5. The second author is supported by S\~ao Paulo Research Foundation (FAPESP) grant 19/13040-7. 

\addcontentsline{toc}{chapter}{Bibliografia}
\bibliographystyle{siam}
\bibliography{Referencias2.bib}
\end{document}